\documentclass[12pt]{article}
\usepackage{amsmath,amsthm,amsfonts,amssymb,amscd,enumerate}
\RequirePackage{graphics}

 \theoremstyle{plain} \newtheorem{theorem}{Theorem}[section]
\newtheorem{proposition}[theorem]{Proposition}

\newtheorem{lemma}[theorem]{Lemma}
\newtheorem{observation}[theorem]{Observation}
\newtheorem*{thmA}{Theorem A} \newtheorem*{thmB}{Theorem B}
 \theoremstyle{definition}
\newtheorem{definition}[theorem]{Definition} \theoremstyle{remark}
\newtheorem{remark}[theorem]{Remark}

\newcommand{\eltwo}{\ell^2} 
\newcommand{\ca}{\mathcal {A}} 
 
 \newcommand{\ch}{\mathcal {H}}
 
\newcommand{\cn}{\mathcal {N}} \newcommand{\ct}{\mathcal {T}}
\newcommand{\cu}{\mathcal {U}} \newcommand{\cv}{\mathcal {V}}
 
\newcommand{\cc}{{\mathbb C}} 
 
\newcommand{\rr}{{\mathbb R}}

\newcommand{\fh}{\mathfrak {H}} 
\newcommand{\sg}{\mathrm {sing}}

\newcommand{\ol}{\overline} \newcommand{\wt}{\widetilde}
\newcommand{\wh}{\widehat}
\newcommand{\ga}{\alpha} \newcommand{\gb}{\beta}

\newcommand{\gs}{\sigma} 
 
\newcommand{\gl}{\lambda}

 \newcommand{\gL}{\Lambda}
\newcommand{\gS}{\Sigma}

\newcommand{\M}{\wt{M}(\ca)} 
 \newcommand{\wu}{\widehat{\cu}}

\newcommand{\Card}{\operatorname{Card}}
\newcommand{\Hom}{\operatorname{Hom}}

\newcommand{\Min}{\operatorname{Min}}
\newcommand{\MU}{\Min (U)}
\newcommand{\rk}{\operatorname{rk}}

\newenvironment{enumeratei}{\begin{enumerate}[\upshape (i)]}%
		{\end{enumerate}}

\begin{document}

\title{The $\eltwo$-cohomology of hyperplane complements}

\author{M. W. Davis \thanks{The first author was partially supported by
NSF grant DMS 0405825.}  \and
T. Januszkiewicz\thanks{The second author was partially supported by
NSF grant DMS 0405825.}  \and I. J. Leary \thanks{The third author was
partially supported by NSF grant DMS 0505471.}  }


\date{\today} \maketitle

\begin{abstract}
We compute the $\eltwo$-Betti numbers of the complement of a 
finite collection of affine hyperplanes in $\cc^n$.  At most 
one of the $\eltwo$-Betti numbers is non-zero.  

\noindent
\textbf{AMS classification numbers}.  Primary: 52B30
Secondary: 32S22, 52C35, 57N65, 58J22.  
\smallskip

\end{abstract}

\section{Introduction}\label{intro}
Suppose $X$ is a finite CW complex with universal cover $\wt{X}$.  For
each $p\ge 0$, one can associate to $X$ a Hilbert space,
$\ch^p(\wt{X})$, the $p$-dimensional ``reduced $\eltwo$-cohomology,'' 
cf.~\cite{eckmann}.  
Each $\ch^p(\wt{X})$ is a unitary $\pi_1(X)$-module.  Using the
$\pi_1(X)$-action, one can attach a nonnegative real number
called ``von Neumann dimension'' to such a Hilbert space.  The
``dimension'' of $\ch^p(\wt{X})$ is called the $p^{\mathrm{th}}$
$\eltwo$-\emph{Betti number of $X$}.

Here we are interested in the case where $X$ is the complement of a
finite number of affine hyperplanes in $\cc^n$.  (Technically, in
order to be in compliance with the first paragraph, we should replace
the complement by a homotopy equivalent finite CW complex.  However, to keep from pointlessly complicating the notation, we shall ignore this technicality.)  Let $\ca$
be the finite collection of hyperplanes, $\gS(\ca)$ their union and
$M(\ca):=\cc^n-\gS(\ca)$.  The \emph{rank} of $\ca$ is the maximum
codimension $l$ of any nonempty intersection of hyperplanes in $\ca$.
It turns out that the ordinary (reduced) homology of $\gS(\ca)$
vanishes except in dimension $l-1$ (cf. Proposition~\ref{p:e}).  Let
$\gb(\ca)$ denote the rank of $\ol{H}_{l-1}(\gS(\ca))$.  Our main
result, proved as Theorem~\ref{t:main}, is the following.

\begin{thmA}
Suppose $\ca$ is an affine hyperplane arrangement of rank $l$.  Only
the $l^{\mathrm{th}}$ $\eltwo$-Betti number of $M(\ca)$ can be nonzero
and it is equal to $\gb(\ca)$.
\end{thmA}

This is reminiscent of a well-known result about the cohomology of
$M(\ca)$ with coefficients in a generic flat line bundle ( ``generic''
is defined in Section~\ref{s:generic}).  This result is proved as
Theorem~\ref{t:generic}.  We state it below.

\begin{thmB}
Suppose that $L$ is a generic flat line bundle over $M(\ca)$.  Then
$H^*(M(\ca);L)$ vanishes except in dimension $l$ and $\dim_\cc
H^l(M(\ca);L)=\gb(\ca)$.
\end{thmB}

Both theorems have similar proofs.  In the case of Theorem A the basic
fact is that the $\eltwo$-Betti numbers of $S^1$ vanish.  (In other
words, if the universal cover $\rr$ of $S^1$ is given its usual cell
structure, then $\ch^*(\rr)=0$.)  Similarly, for Theorem B, if $L$ is
a flat line bundle over $S^1$ corresponding to an element $\gl\in
\cc^*$, with $\gl\neq 1$, then $H^*(S^1;L)=0$.  By the K\"unneth
Formula, there are similar vanishing results for any central
arrangement.  To prove the general results, one considers an open
cover of $M(\ca)$ by ``small'' open sets each homeomorphic to
the complement of a central arrangement.  The $E_1$-page of the
resulting Mayer-Vietoris spectral sequence is nonzero only along the
bottom row, 
where it can be identified with the simplicial cochains with constant
coefficients on a pair $(N(\cu),N(\cu_\sg))$, which is homotopy equivalent
to $(\cc^n,\gS)$.  It follows that the $E_2$-page can be nonzero only
in position $(l,0)$.  (Actually, in the case of Theorem A, technical
modifications must be made to the above argument.  Instead of reduced
$\eltwo$-cohomology one takes local coefficients in the von Neumann
algebra associated to the fundamental group and the vanishing results
only hold modulo modules which do not contribute to the $\eltwo$-Betti
numbers.)

In \cite{dl} the first and third authors proved a similar result for
the $\eltwo$-cohomology of the universal cover of the Salvetti complex
associated to an arbitrary Artin group (as well as a formula for the
cohomology of the Salvetti 
complex with generic,
$1$-dimensional local coefficients). This can be interpreted as a
computation of the $\eltwo$-cohomology of universal covers of
hyperplane complements associated to infinite reflection groups.
Although the main argument in \cite{dl} uses an explicit description
of the chain complex of the Salvetti complex, an alternative argument,
similar to the one outlined above, is given in \cite[Section 10]{dl}.

We thank the referee for finding some mistakes in the first version of this paper.

\section{Hyperplane arrangements}\label{s:hyper}
A \emph{hyperplane arrangement} $\ca$ is a finite collection of affine
hyperplanes in $\cc^n$.  A \emph{subspace} of $\ca$ is a nonempty
intersection of hyperplanes in $\ca$.  Denote by $L(\ca)$ the poset of
subspaces, partially ordered by inclusion, and let $\ol{L}(\ca):=L(\ca)\cup \{\cc^n\}$.  An arrangement is
\emph{central} if $L(\ca)$ has a minimum element.  Given $G\in
L(\ca)$, its \emph{rank}, $\rk(G)$, is the codimension of $G$ in
$\cc^n$.  The minimal elements of $L(\ca)$ are a family of parallel
subspaces and they all have the same rank.  The \emph{rank} of an
arrangement $\ca$ is the rank of a minimal element in $L(\ca)$.  $\ca$
is \emph{essential} if $\rk(\ca)=n$.

The \emph{singular set} $\gS(\ca)$ of the arrangement is the union of
hyperplanes in $\ca$ (so that $\gS(\ca)$ is a subset of $\cc^n$).  The
complement of $\gS(\ca)$ in $\cc^n$ is denoted $M(\ca)$.  When there
is no ambiguity we will drop the ``$\ca$'' from our notation and write
$L$, $\gS$ or $M$ instead of $L(\ca)$, $\gS(\ca)$ or $M(\ca)$.  

\begin{proposition}\label{p:e}
$\gS$ is homotopy equivalent to a wedge of $(l-1)$-spheres,
where $l=\rk(\ca)$.  (So, if $\ca$ is essential, the spheres are
$(n-1)$-dimensional.)
\end{proposition}

\begin{proof}
The proof follows from the usual ``deletion-restriction''
argument and induction.  If the rank $l$ is $1$, 
then $\gS$ is the disjoint union of a finite family of parallel hyperplanes.  Hence, $\gS$ is homotopy equivalent to a finite set of points, i.e., to a wedge of $0$-spheres.  Similarly, when
$l=2$, it is easy to see that $\gS$ is homotopy equivalent to a connected graph; hence, a wedge of
$1$-spheres.  So, assume by induction that $l>2$.  Choose a hyperplane
$H\in \ca$, let $\ca'=\ca-\{H\}$ and let $\ca''$ be the restriction of
$\ca$ to $H$ (i.e., $\ca'':=\{H'\cap H\mid H'\in \ca'\}$).  Put
$\gS'=\gS(\ca')$, $\gS''=\gS(\ca'')$, $l'=\rk(\ca')$ and
$l''=\rk(\ca'')$. We can also assume by induction on $\Card(\ca)$ that
$\gS'$ and $\gS''$ are homotopy equivalent to wedges of spheres. If
$l'<n$ and $H$ is transverse to the minimal elements of $L(\ca')$, then
$l''=l$, the arrangement splits as a product, $\gS=\gS''\times \cc$,
and we are done by induction.  In all other cases $l'=l$ and
$l''=l-1$.  We have $\gS=\gS' \cup H$ and $\gS'\cap H=\gS''$.  $H$ is
simply connected and since $l>2$, $\gS'$ is simply connected and
$\gS''$ is connected.  By van Kampen's Theorem, $\gS$ is simply
connected.  Consider the exact sequence of the pair $(\gS,\gS')$:
\[
\to H_*(\gS')\to H_*(\gS)\to H_*(\gS,\gS') \to .
\]
There is an excision isomorphism, $H_*(\gS,\gS')\cong H_*(H,\gS'')$.  
Since $H$ is contractible it follows that 
$H_*(H,\gS'')\cong \ol{H}_{*-1}(\gS'')$.  By induction,
$\ol{H}_*(\gS')$ is concentrated in dimension $l-1$ and
$\ol{H}_*(\gS'')$ in dimension $l-2$.  So, $\ol{H}_*(\gS)$ is also
concentrated in dimension $l-1$.  It follows that $\gS$ is homotopy
equivalent to a wedge of $l-1$ spheres.
\end{proof}

\section{Certain covers and their nerves}\label{s:covers}
Suppose $\cu=\{U_i\}_{i\in I}$ is a cover of some space $X$ (where $I$
is some index set).  Given a subset $\gs\subset I$, put
$U_\gs:=\bigcap_{i\in \gs} U_i$.  Recall that the \emph{nerve} of
$\cu$ is the simplicial complex $N(\cu)$, defined as follows.  Its
vertex set is $I$ and a finite, nonempty subset $\gs\subset I$ spans a
simplex of $N(\cu)$ if and only if $U_\gs$ is nonempty.

We shall need to use the following well-known lemma several times in
the sequel, see~\cite[Cor.\ 4G.3 and Ex.\ 4G(4)]{hatcher}

\begin{lemma}\label{l:nerve}
Let $\cu$ be a cover of a paracompact space $X$ and suppose that
either (a) each $U_i$ is open or (b) $X$ is a CW complex and each
$U_i$ is a subcomplex.  Further suppose that for each simplex $\gs$ of
$N(\cu)$, $U_\gs$ is contractible.  Then $X$ and $N(\cu)$ are homotopy
equivalent.
\end{lemma}

Suppose $\ca$ is a hyperplane arrangement in $\cc^n$.
An open convex subset $U$ in $\cc^n$ is \emph{small} (with respect to $\ca$) if the following two conditions hold:
\begin{enumeratei}
\item
$\{G\in \ol{L}(\ca)\mid G\cap U\neq\emptyset\}$ has a unique minimum element $\MU$.
\item
A hyperplane $H\in\ca$ has nonempty intersection with $U$ if and only if $\MU \subset H$.
\end{enumeratei}
The intersection of two small convex open sets is also small; hence, the same is true for any finite intersection of such sets.

Now let $\cu=\{U_i\}_{i\in I}$ be an open cover of $\cc^n$ by small convex sets.  Put
\[
\cu_\sg:= \{U\in \cu\mid U\cap \gS\neq \emptyset\}.
\]

\begin{lemma}\label{l:nu}
$N(\cu)$ is a contractible simplicial complex and $N(\cu_\sg)$ is a
subcomplex homotopy equivalent to $\gS$.  Moreover, $H_*(N(\cu),N(\cu_\sg))$ is concentrated in
dimension $l$, where $l=\rk \ca$.
\end{lemma}

\begin{proof}
$\cu_\sg$ is an open cover of a neighborhood of $\gS$ which deformation retracts onto $\gS$.  For each simplex $\gs$ of $N(\cu)$, $U_\gs$ is contractible (in fact, it is a small convex open set).
By Lemma~\ref{l:nerve}, $N(\cu)$ is homotopy equivalent to $\cc^n$ and $N(\cu_\sg)$ is homotopy equivalent to $\gS$.  The last sentence of the lemma  follows from Proposition~\ref{p:e}.
\end{proof}

\begin{remark}
Lemma~\ref{l:nerve} can also be used to show that the geometric realization of $L$ is homotopy equivalent to $\gS$.
\end{remark}

\begin{definition}\label{d:beta}
$\gb(\ca)$ is the rank of $H_l(N(\cu), N(\cu_\sg))$.
\end{definition}

Equivalently, $\gb(\ca)$ is the rank of $H_l(\cc^n,\gS(\ca))$ (or of
$\ol{H}_{l-1}(\gS(\ca))$.  Also, it is not difficult to see that
$(-1)^l\gb(\ca)=\chi(\cc^n,\gS)=1-\chi(\gS)=\chi(M)$, where $\chi(\,
)$ denotes the Euler characteristic.

\begin{remark}
Suppose $\ca_\rr$ is an arrangement of real hyperplanes in $\rr^n$ and
$\gS_\rr\subset \rr^n$ is the singular set.  Then $\rr^n-\gS_\rr$ is a
union of open convex sets called \emph{chambers} and $\gb(\ca_\rr)$ is
the number of bounded chambers.  If $\ca$ is the complexification of
$\ca_\rr$, then $\gS(\ca)\sim\gS(\ca_\rr)$.  Hence,
$\gb(\ca)=\gb(\ca_\rr)$.
\end{remark}

For any small open convex set $U$, put
\[
\wh{U}:=U-\gS(\ca) = U\cap M(\ca).
\]
Since $U$ is convex, $(U,U\cap\,\gS(\ca))$ is homeomorphic to $(\cc^n,\gS(\ca_G))$, where $G=\MU$ and $\ca_G$ is the central subarrangement defined by
\[
\ca_G:=\{H\in \ca\mid G\subset H\}.  
\]
($G$ might be $\cc^n$, in which case $\ca_G=\emptyset$.) Hence, $\wh {U}$ is homeomorphic to $M(\ca_G)$, the complement of a central subarrangement.

The next lemma is well-known.

\begin{lemma}\label{l:small}
Suppose $U$ is a small open convex set. Then $\pi_1(\wh{U})$ is a retract of
$\pi_1(M(\ca))$.
\end{lemma}

\begin{proof}
The composition of the two inclusions, $\wh{U}\hookrightarrow
M(\ca)\hookrightarrow M(\ca_G)$ is a homotopy equivalence, where $G=\MU \in L(\ca)$.
\end{proof}

By intersecting the elements of $\cu$ with $M$ ($=\cc^n-\gS$) we get
an induced cover $\wu$ of $M$.  An element of $\wu$ is a deleted small
convex open set $\wh{U}$ for some $U\in \cu$.  Similarly, by
intersecting $\cu_\sg$ with $M$ we get an induced cover $\wu_\sg$ of
a deleted neighborhood of $\gS$.  The key observation is the following.

\begin{observation}\label{o}
$N(\wh{\cu})=N(\cu)$ and $N(\wh{\cu}_\sg)=N(\cu_\sg)$.
\end{observation}

\section{The Mayer-Vietoris spectral sequence}\label{s:mv}
Let $X$ be a space, $\pi=\pi_1(X)$ and $r:\wt{X}\to X$ the universal
cover.  Given a left $\pi$-module $A$, define
\[
C^*(X;A):=\Hom_\pi(C_*(\wt{X}),A),
\]
the cochains with \emph{local coefficients in $A$}.  Taking cohomology
gives $H^*(X;A)$.

Let $\cu$ be an open cover of $X$ and $N=N(\cu)$ its nerve.  Let
$N^{(p)}$ denote the set of $p$-simplices in $N$.  There is an induced
cover $\wt{\cu}:=\{r^{-1}(U)\}_{U\in \cu}$ with the same
nerve.  There is a Mayer-Vietoris double complex
\[
C_{p,q}=\bigoplus_{\gs\in N^{(p)}}C_q(r^{-1}(U_\gs))
\]
(cf. \cite[\S VII.4]{brown}) and a corresponding double cochain
complex with local coefficients:
\[
C^{p,q}(A):=\Hom_\pi(C_{p,q};A).
\]
The cohomology of the total complex is $H^*(X;A)$.  Now suppose that for each simplex $\gs$ of $N$, $U_\gs$ is connected and that $\pi_1(U_\gs)\to
\pi_1(X)$ is injective. (This implies that $r^{-1}(U_\gs)$ is a disjoint union of copies of the universal cover $\wt{U}_\gs$.)
We get a
spectral sequence with $E_1$-page
\begin{equation}\label{e:E1}
E^{p,q}_1=\bigoplus_{\gs\in N^{(p)}}H^q(U_\gs;A).
\end{equation}
Here $H^q(U_\gs;A)$ means the cohomology of
$\Hom_\pi(C_*(r^{-1}(U_\gs)),A)$ or equivalently,  of
$\Hom_{\pi_1(U_\gs)} (C_*(\wt{U}_\gs); A)$.  The $E_2$-page has the form $E^{p,q}=H^p(N;\fh^q)$, where $\fh^q$ means the functor $\gs\to H^q(U_\gs;A)$.  The spectral 
sequence converges to $H^*(X;A)$.

In the next two sections we will apply this spectral sequence to the
case where $X$ is $M(\ca)$ and the open cover is $\wu$ from the previous
section.  By Lemma~\ref{l:small}, $\pi_1(\wh{U}_\gs) \to \pi_1(M(\ca)$ is injective so we get a spectral sequence with $E_1$-page given by (\ref{e:E1}).
Moreover, the $\pi$-module $A$ will be such that for any
simplex $\gs$ in $N(\wh{\cu}_\sg)$, $H^q(\wh{U}_\gs;A)=0$ for all $q$
(even for $q=0$) while for a simplex $\gs$ of $N(\wh{\cu})$ which is
not in $N(\wh{\cu}_\sg)$, $H^q(U_\gs;A)=0$ for all $q>0$ and is constant
(i.e., independent of $\gs$) for $q=0$.  Thus $E^{p,q}_1$ will vanish
for $q>0$ and $E^{*,0}_1$ can be identified with the cochain complex
$C^*(N(\cu), N(\cu_\sg))$ with constant coefficients.

\section{Generic coefficients}\label{s:generic}
Here we will deal with $1$-dimensional local coefficient systems.  We
begin by considering such local coefficients on $S^1$.  Let $\ga$ be a
generator of the infinite cyclic group $\pi_1(S^1)$.  Suppose $k$ is a
field of characteristic $0$ and $\gl\in k^*$.  Let $A_\gl$ be the
$k[\pi_1(S^1)]$-module which is a $1$-dimensional $k$-vector space on
which $\ga$ acts by multiplication by $\gl$.

\begin{lemma}\label{l:s1}
If $\gl\neq 1$, then $H^*(S^1;A_\gl)$ vanishes identically.
\end{lemma}

\begin{proof}
If $S^1$ has its usual CW structure with one $0$-cell and one
$1$-cell, then in the chain complex for its universal cover both $C_0$
and $C_1$ are identified with the group ring $k[\pi_1(S^1)]$ and the
boundary map with multiplication by $1-t$, where $t$ is the generator
of $\pi_1(S^1)$.  Hence, the coboundary map $C^0(S^1;A_\gl)\to
C^1(S^1;A_\gl)$ is multiplication by $1-\gl$.  
\end{proof}

Next, consider $M(\ca)$.  Its fundamental group $\pi$ is generated by
loops $a_H$ for $H\in \ca$, where the loop $a_H$ goes once around the 
hyperplane $H$ in the ``positive'' direction.  
Let $\ga_H$ denote the image of $a_H$ in
$H_1(M(\ca))$.  Then $H_1(M(\ca))$ is free abelian with basis
$\{\ga_H\}_{H\in \ca}$.  So, a homomorphism $H_1(M(\ca))\to k^*$ is
determined by an $\ca$-tuple $\gL\in (k^*)^\ca$, where
$\gL=(\gl_H)_{H\in \ca}$ corresponds to the homomorphism sending
$\ga_H$ to $\gl_H$.  Let $\psi_\gL:\pi\to k^*$ be the composition of
this homomorphism with the abelianization map $\pi\to H_1(M(\ca))$.
The resulting local coefficient system on $M(\ca)$ is denoted $A_\gL$.
The next lemma follows from Lemma~\ref{l:s1}.

\begin{lemma}\label{l:central}
Suppose $\ca$ is a nonempty central arrangement and $\gL$ is such that
$\prod_{H\in\ca} \gl_H \neq 1$.  Then $H^q(M(\ca))$ vanishes for all
$q$.
\end{lemma}

\begin{proof}
Without loss of generality we can suppose the elements of $\ca$ are
linear hyperplanes.  The Hopf bundle $M(\ca)\to M(\ca)/S^1$ is trivial
(cf. \cite[Prop. 5.1, p. 158]{ot}); so, $M(\ca)\cong B\times S^1$,
where $B=M(\ca)/S^1$.  Let $i:S^1\to M(\ca)$ be inclusion of the
fiber. The induced map on $H_1(\ )$ sends $\ga$ to $\sum \ga_H$.
Thus, if we pull back $A_\gL$ to $S^1$, we get $A_\gl$, where
$\gl=\prod_{H\in\ca} \gl_H$.  The condition on $\gL$ is $\gl\neq 1$,
which by Lemma~\ref{l:s1} implies that $H^*(S^1;A_\gl)$ vanishes
identically.  By the K\"unneth Formula $H^*(M(\ca);A_\gL)$ also
vanishes identically.
\end{proof}

Returning to the case where $\ca$ is a general arrangement, for each
simplex $\gs$ in $N(\wu)$, let $\ca_\gs:=\ca_{\Min(U_\gs)}$ be the corresponding central
arrangement (so that $\wh{U}_\gs \cong M(\ca_\gs)$).  Given $\gL\in
(k^*)^\ca$, put
\[
\gl_\gs:=\prod_{H\in \ca_\gs} \gl_H.
\]
Call $\gL$ \emph{generic} if $\gl_\gs\neq 1$ for all $\gs\in N(\cu_\sg)$.

\begin{theorem}\label{t:generic} \textup{(Compare \cite[Thm. 4.6,
      p. 160]{vs}).} 
Let $\ca$ be an affine arrangement of rank $l$ and $\gL$ a generic
$\ca$-tuple in $k^*$.  Then $H^*(M(\ca);A_\gL)$ is concentrated in
degree $l$ and $$\dim_k H^l(M(\ca);A_\gL)=\gb(\ca).$$
\end{theorem}

\begin{proof}
We have an open cover of $\M$, $\{r^{-1}(\wh{U})\}_{U\in \cu}$.  By Observation~\ref{o}, its nerve
is $N(\cu)$.  By Lemma~\ref{l:central} and the last paragraph of
Section~\ref{s:mv}, the $E_1$-page of the Mayer-Vietoris spectral
sequence is concentrated along the bottom row where it can be
identified with $C^*(N(\cu),N(\cu_\sg);k)$.  So, the $E_2$-page is concentrated
on the bottom row and $E^{p,0}_2=H^p(N(\cu),N(\cu_\sg);k)$.  By
Lemma~\ref{l:nu}, this group is nonzero only for $p=l$ and 
\[
\dim_k E^{l,0}_2=\dim_k H^l(N(\cu), N(\cu_\sg);k)=\gb(\ca).
\]
\end{proof}

\begin{remark}
When $k=\cc$, a $1$-dimensional local coefficient system on~$X$ 
is the same thing as a flat line bundle over $X$.
\end{remark}

\section{$\eltwo$-cohomology}\label{s:ltwo}
For a discrete group $\pi$, $\eltwo\pi$ denotes the Hilbert space of
complex-valued, square integrable functions on $\pi$.  There are
unitary $\pi$-actions on $\eltwo\pi$ by either left or right
multiplication; hence, $\cc\pi$ acts either from the left or right as
an algebra of operators.  The \emph{associated von Neumann algebra}
$\cn\pi$ is the commutant of $\cc \pi$ (acting from, say, the right on
$\eltwo\pi$).

Given a finite CW complex $X$ with fundamental group $\pi$, the space
of $\eltwo$-cochains on its universal cover $\wt{X}$ is the same as
$C^*(X;\eltwo\pi)$, the cochains with local coefficients in $\eltwo
\pi$.  The image of the coboundary map need not be closed; hence,
$H^*(X;\eltwo\pi)$ need not be a Hilbert space.  To remedy this, one
defines the \emph{reduced} $\eltwo$-cohomology $\ch^*(\wt{X})$ to be
the quotient of the space of cocycles by the closure of the space of
coboundaries.  We shall also use the notation $\ch^*(X;\eltwo\pi)$ for
the same space.

The von Neumann algebra admits a trace.  Using this, one can attach a
``dimension,'' $\dim_{\cn\pi} V$, to any closed, $\pi$-stable subspace
$V$ of a finite direct sum of copies of $\eltwo\pi$ (it is the trace
of orthogonal projection onto $V$).  The nonnegative real number
$\dim_{\cn\pi}(\ch^p(X;\eltwo\pi))$ is the $p^{\mathrm {th}}$
\emph{$\eltwo$-Betti number} of~$X$.

A technical advance of L\"uck \cite[Ch.~6]{luckbk} is the use local
coefficients in $\cn\pi$ in place of the previous version of
$\eltwo$-cohomology.  
He shows there is a well-defined
dimension function on $\cn\pi$-modules, $A\to \dim_{\cn\pi} A$, which gives the same gives the same answer
for $\eltwo$-Betti numbers, i.e., for each~$p$ 
one has that 
$\dim_{\cn\pi}H^p(X;\cn\pi)=\dim_{\cn\pi}\ch^p(X;\eltwo\pi)$.  
Let $\ct$ be the class of $\cn \pi$-modules of dimension $0$.  
The dimension function is additive with respect to short exact sequences.  
This allows one to define $\eltwo$-Betti numbers 
for spaces more general than finite complexes.
The class $\ct$ 
is a Serre class of $\cn\pi$-modules~\cite{jps}, which allows one
to compute $\eltwo$-Betti numbers by working with spectral 
sequences modulo~$\ct$.


\begin{lemma}\label{l:central2}
Suppose $\ca$ is a nonempty central arrangement.  Then, for all $q\ge
0$, $H^q(M(\ca);\cn\pi)$ lies in $\ct$.  In other words, all
$\eltwo$-Betti numbers of $M(\ca)$ are zero.
\end{lemma}

\begin{proof}
The proof is along the same line as that of Lemma~\ref{l:central}.  It
is well-known that the reduced $\eltwo$-cohomology of $\rr$ vanishes. 
Since $M(\ca)=S^1\times B$, the result follows from the K\"unneth
Formula for $\eltwo$-cohomology in~\cite[6.54~(5)]{luckbk}.  
\end{proof}

\begin{theorem}\label{t:main}
Suppose $\ca$ is an affine hyperplane arrangement.  Then
\[
H^*(M(\ca);\cn\pi)\cong H^*(N(\cu), N(\cu_\sg))\otimes \cn\pi \pmod {\ct}
\]
Hence, for $l=\rk(\ca)$, the $\eltwo$-Betti numbers of $M(\ca)$
vanish except in dimension $l$, where
$\dim_{\cn\pi}\ch^l(\wt{M}(\ca))=\gb(\ca)$.
\end{theorem}

\begin{proof}
For each $\gs\in N(\cu_\sg)$, let $\pi_\gs:=\pi_1(U_\gs)$.  By Lemma~\ref{l:central2}, 
\[
\dim_{\cn \pi_\gs} H^*(M(\ca_\gs); \cn\pi_\gs)=0.
\]
Since the $\cn_\pi$-module $H^*(M(\ca_\gs), \cn\pi)$ is induced from $H^*(M(\ca_\gs), \cn\pi)$,
\[
\dim_{\cn\pi}H^*(M(\ca_\gs); \cn\pi)=\dim_{\cn\pi_\gs}H^*(M(\ca_\gs); \cn\pi_\gs)
=0.
\]
As in the proof of Theorem~\ref{t:generic}, it follows that the $E_1$-page of the spectral 
sequence consists of modules in $\ct$, except that $E_1^{*,0}$ 
is identified with $C^*(N(\cu),N(\cu_\sg))\otimes \cn(\pi)$.  Similarly, 
the $E_2$-page consists of modules in $\ct$, except that 
$E_2^{*,0}$ is identified with $H^*(N(\cu),N(\cu_\sg))\otimes \cn\pi$.  
For each subsequent differential, either the source or the target
is a module in $\ct$, and hence for each $i$ and $j$ one has that 
$E_\infty^{i,j} \cong E_2^{i,j} \pmod \ct$.  
The claim follows since the filtration of $H^*(M(\ca);\cn\pi)$ given 
by the $E_\infty$-page of the spectral sequence is finite.  
\end{proof}


\obeylines
M. W. Davis {\tt mdavis@math.ohio-state.edu }
T. Januszkiewicz {\tt tjan@math.ohio-state.edu} 
I. J. Leary {\tt leary@math.ohio-state.edu} 
\

Department of Mathematics, 
The Ohio State University, 
231 W. 18th Ave., 
Columbus 
Ohio 43210

\end{document}